\documentclass[12pt,a4paper]{amsart}
\usepackage{amssymb,amsfonts,amsthm,amsmath}
\usepackage{graphicx} % for inclusion of image
\usepackage{color}

\newcommand{\summ}{\mathop{\sum}}

\theoremstyle{plain}

\numberwithin{equation}{section} \numberwithin{theorem}{section}
\numberwithin{lemma}{section} \numberwithin{definition}{section}
\numberwithin{corollary}{section} \textheight =24cm
\begin{document}
\title[Ramanujan's nested radicals]{Variations on Ramanujan's nested radicals}
\author{Geoffrey B Campbell}
\address{Mathematical Sciences Institute \\
         The Australian National University \\
         ACT, 0200, Australia}
 \email{Geoffrey.Campbell@anu.edu.au}

\author{Aleksander Zujev}
\address{Physics Department \\
         University of California \\
         Davis, California 95616, USA}
\email{azujev@ucdavis.edu}

\begin{abstract}
We give new nested radical equations of similar kind to Ramanujan's questions to the Indian Mathematical Society 100 years ago. While many have since considered these from the perspectives of the Notebooks of Ramanujan and from the theory of Class numbers and Units, there seems no comprehensive theory to cover off the results, and it seems always possible to find new and surprizing elementary equations of this kind. We also consider a few more methods of constructing nested radicals.
%%Finally, we consider an example of infinitely nested radicals.
\end{abstract}
\maketitle

\section{Introduction} \label{S:intro}
Arguably the most cited examples of Ramanujan's nested radicals are
\begin{equation} \label{E:1.1}
              \sqrt[4]{\frac{3+2\sqrt[4]{5}}{3-2\sqrt[4]{5}}} = \frac{\sqrt[4]{5} +1}{\sqrt[4]{5}-1}
\end{equation}
and
\begin{equation} \label{E:1.2}
              \sqrt[3]{\sqrt[3]{2}-1} = \sqrt[3]{\frac{1}{9}} - \sqrt[3]{\frac{2}{9}} + \sqrt[3]{\frac{4}{9}}.
\end{equation}

%% The right side of (\ref{E:1.2}) is also equal to:
%% \begin{eqnarray} \label{E:1.3}&& \sqrt[24]{1+100\sqrt[3]{2}-80\sqrt[3]{4}}; \\
%% \label{E:1.4} && \sqrt[27]{180\sqrt[3]{4}-99\sqrt[3]{2}-161}; \\
%% \label{E:1.5} && \sqrt[36]{1513-1662\sqrt[3]{2}-366\sqrt[3]{4}}; \\
%% \label{E:1.6} && \sqrt[54]{-45359+96678\sqrt[3]{2}-48159\sqrt[3]{4}}; \\
%% \label{E:1.7} && \sqrt[81]{51642361-28411857\sqrt[3]{2}-9982143\sqrt[3]{4}}.
%% \end{eqnarray}

In particular, (\ref{E:1.2}) has been used and reused in texts as an exercise problem, and used and reused in Mathematical Olympiad questions at various times. In this paper we first look at equations (\ref{E:1.1}) and (\ref{E:1.2}) and their possible extensions. Then we give examples of new nested identities of this kind, and present a method to construct such results. Along the way we try to examine what exactly makes such identities "interesting".

\section{Ramanujan-style nested radicals}

\subsection{Equation (\ref{E:1.1}).}

Let's consider a general equation in the shape of equation (\ref{E:1.1}):
\begin{equation} \label{E:2.1}
              \sqrt[4]{\frac{x+y\sqrt[4]{b}}{x-y\sqrt[4]{b}}} = \frac{z+w\sqrt[4]{b}}{z-w\sqrt[4]{b}}.
\end{equation}
Solving it, we come to a Diophantine equation
\begin{equation} bw^4 \; = \; 5z^4,  \end{equation}
with its only solution $b = 5$.
So Ramanujan's equation (\ref{E:1.1}) is the only one possible for radicals of the form (\ref{E:2.1}).

A similar equation with different radical form,
\begin{equation} \label{E:2.1a}
              \sqrt[n]{\frac{x+y\sqrt[n]{b}}{x-y\sqrt[n]{b}}} = \frac{z+w\sqrt[n]{b}}{z-w\sqrt[n]{b}}
\end{equation}
also doesn't give solutions.

The even more general equation
\begin{equation} \label{E:2.1b}
              \sqrt[n]{\frac{x+y\sqrt[m]{b}}{x-y\sqrt[m]{b}}} = \frac{z+w\sqrt[m]{b}}{z-w\sqrt[m]{b}}
\end{equation}
has some solutions, for example
\begin{eqnarray} \label{E:2.1c}
              \sqrt[4]{\frac{7+4\sqrt[]{3}}{7-4\sqrt[]{3}}} = \frac{3+\sqrt[]{3}}{3-\sqrt[]{3}}
\end{eqnarray}

\subsection{Equation (\ref{E:1.2}).}

If we closely consider RHS of the equation (\ref{E:1.2}), we can see that it's three terms of the geometric series.
It is easy to construct extensions of the equation (\ref{E:1.2}), taking different sub-sequences of the series:

\begin{eqnarray} \label{E:2.2}
\sqrt[3]{\frac{1}{16} (-7+3 \sqrt[3]{2} + 6 \sqrt[3]{4})} = 1 - \sqrt[3]{\frac{1}{2}} +\sqrt[3]{\frac{1}{4}}
-\sqrt[3]{\frac{1}{8}} + \sqrt[3]{\frac{1}{16}},
\end{eqnarray}
\begin{eqnarray} \label{E:2.3}
\sqrt[3]{\frac{9}{32} (\sqrt[3]{2} -1)} = 1 - \sqrt[3]{\frac{1}{2}} +\sqrt[3]{\frac{1}{4}} -\sqrt[3]{\frac{1}{8}}
+\sqrt[3]{\frac{1}{16}} -\sqrt[3]{\frac{1}{32}}.
\end{eqnarray}

%% (1/16 (-7+3 2^(1/3)+6 2^(2/3)))^(1/3) = 1^(1/3) - (1/2)^(1/3) +(1/4)^(1/3) - (1/8)^(1/3) + (1/16)^(1/3)
%% (9/32 (2^(1/3)-1))^(1/3) = 1^(1/3) - (1/2)^(1/3) + (1/4)^(1/3) - (1/8)^(1/3) + (1/16)^(1/3) - (1/32)^(1/3)
%% (9 (2^(1/3)-1))^(1/3) = Sum_{k=0}^2 1 - 2^(1/3) + 4^(1/3)

They are all cases of a geometric series sum
\begin{eqnarray} \label{E:2.4}
\sqrt[3]{\frac{1}{27}(-(-2)^m+1)^3 (\sqrt[3]{2}-1)} = \summ_{k=0}^{3m-1} \frac{(-2)^{\frac{k}{3}}}{9^{\frac{1}{3}}}.
\end{eqnarray}

%% ((-(-2)^(m)+1)^3 (2^(1/3)-1))^(1/3) = Sum_{k=0}^{3m-1} (-2)^(k/3)/9^(1/3)
Ramanujan's formula (\ref{E:1.2})
%% (2^(1/3)-1)^(1/3) = Sum_{k=0}^2 (1/9)^(1/3) - (2/9)^(1/3) + (4/9)^(1/3)
is just the first three terms of the series.

Going the other direction into negative exponents, we get
\begin{eqnarray} \label{E:2.5}
\sqrt[3]{\frac{(2^m+(-1)^{m+1})^3}{27 \cdot 2^{3m-1}} (\sqrt[3]{2}-1)} = \summ_{k=-(3m-1)}^0 \frac{(-2)^{\frac{k}{3}}}{9^\frac{1}{3}},
\end{eqnarray}

and
\begin{eqnarray} \label{E:2.6}
\sqrt[3]{\frac{2}{27}(\sqrt[3]{2}-1)} = \summ_{k=-\infty}^0 \frac{(-2)^{\frac{k}{3}}}{9^\frac{1}{3}}.
\end{eqnarray}

Hence, we can state a few equations similar to (\ref{E:1.2}), that are not quite geometric series:

\begin{eqnarray} \label{E:2.7}
\sqrt[3]{18(1-\sqrt[3]{2})} = \sqrt[3]{2^2} - 2 - \sqrt[3]{2},
\end{eqnarray}
\begin{eqnarray} \label{E:2.8}
\sqrt[3]{-13122(1-\sqrt[3]{2})} = - 9 \sqrt[3]{2^2} + 18 + 9 \sqrt[3]{2},
\end{eqnarray}
\begin{eqnarray} \label{E:2.9}
\sqrt[3]{486(1-\sqrt[3]{2})} = 3 \sqrt[3]{2^2} - 6 - 3 \sqrt[3]{2},
\end{eqnarray}
\begin{eqnarray} \label{E:2.10}
\sqrt[3]{28917 + 64638 \sqrt[3]{7}} = 12 \sqrt[3]{7^2} + 21 - 27 \sqrt[3]{7},
\end{eqnarray}
\begin{eqnarray} \label{E:2.11}
\sqrt[]{406 + 84 \sqrt[4]{7} -90 \sqrt{7}} = \sqrt[4]{7} + 7 \sqrt[]{7} + \sqrt[4]{7^3} - 7.
\end{eqnarray}
These latter two results do not seem to be in the literature, but the other results from (\ref{E:2.5}) to (\ref{E:2.9}) all simplify back to essentially Ramanujan's (\ref{E:1.2}).

\subsection{New nested radicals identities.}

Here are some new Ramanujan-style nested radicals that seem to be different to identities elsewhere:
\begin{eqnarray} \label{eqa}
&& \sqrt[]{1 -  \sqrt[3]{\frac{1}{2}} +5\sqrt[3]{2}+3\sqrt[6]{32}} = 1 +\frac{1}{\sqrt[6]{2}} -\sqrt[6]{2}+\sqrt[3]{2}+\sqrt{2}, \\
\label{eqb} && \sqrt{2-\frac{2}{\sqrt[3]{3}} +\frac{8}{\sqrt[6]{3}}+5\sqrt[3]{3}} =
1+\frac{1}{\sqrt[6]{3}}-\sqrt[6]{3}+\sqrt[3]{3}+\sqrt{3}, \\
\label{eqc} && \sqrt{1+\frac{7}{\sqrt[3]{4}}-\frac{1}{\sqrt[3]{2}+3\sqrt{2}}} =
1+\sqrt[6]{\frac{1}{2}}-\sqrt[6]{2}+\sqrt[3]{\frac{1}{2}}+\sqrt{2}, \\
\label{eqd} && \sqrt{10+8\sqrt[6]{3}-3\sqrt[3]{3}+2\sqrt[3]{9}} =
1+\sqrt[6]{3}+\sqrt[3]{3}-\sqrt{3} +\sqrt[6]{243}.
\end{eqnarray}

A nested radical is "interesting" if there are less terms under the radical left side of the equation,
than you would normally expect from expansion of the terms on the right side.
Ramanujan's example (\ref{E:1.2}) has two terms under the radical on the left and three terms on the right, which on first sight seems surprising, as many terms cancel to enable this.
In (\ref{eqa}) - (\ref{eqd}) there are four terms under the radical on the right and five terms on the left.
Here are two more new "interesting" equations:
\begin{eqnarray}
\label{eqe}&&\sqrt{43+12\sqrt[4]{2}-18\sqrt[4]{8}} = 4\sqrt[4]{2}-3\sqrt[4]{4}+\sqrt[4]{8}+3, \\
\label{eqf}&&\sqrt{7+6\sqrt[6]{2}-\sqrt[3]{2}+\sqrt[3]{4}} = 1 -\sqrt[]{2} +\sqrt[6]{2}+\sqrt[3]{2}+\sqrt[6]{32}.
\end{eqnarray}

Another identity:
\begin{eqnarray}
\sqrt[4]{\sqrt[3]{\frac{1}{9}}+2\sqrt[3]{\frac{2}{9}}-2\sqrt[3]{\frac{4}{9}}} =
\sqrt[5]{\sqrt[3]{9}-\sqrt[3]{\frac{2}{3}}-\sqrt[3]{\frac{4}{3}}} =
\sqrt[6]{1-2\sqrt[3]{2}+\sqrt[3]{4}},
\end{eqnarray}
the leftmost surd being a trivial square under the sixth root.

From Ramanujan's equation (\ref{E:1.2}) we can explore new higher roots by expanding the right side of (\ref{E:1.2}) in any computing program that expands algebraic expressions such as Mathematica, Maple and WolframAlpha. Hence, we discover easily that the LHS of (\ref{E:1.2}) namely, $\sqrt[3]{\frac{1}{9}}-\sqrt[3]{\frac{2}{9}}+\sqrt[3]{\frac{4}{9}}$, may be equivalently written as:
\begin{eqnarray}
\label{eqx2} && \sqrt[8]{4 \sqrt[3]{\frac{2}{3}} - 5 \sqrt[3]{\frac{1}{3}}}, \\
\label{eqx3} && \sqrt[15]{19 - 5\sqrt[3]{2} -8 \sqrt[3]{4}}, \\
\label{eqexa} && \sqrt[24]{1+100\sqrt[3]{2}-80\sqrt[3]{4}}, \\
&& \sqrt[27]{180\sqrt[3]{4}-99\sqrt[3]{2}-161}, \\
&& \sqrt[36]{1513-1662\sqrt[3]{2}+366\sqrt[3]{4}}, \\
&& \sqrt[54]{-45359+96678\sqrt[3]{2}-48159\sqrt[3]{4}}, \\
&& \sqrt[81]{51642361-28411857\sqrt[3]{2}-9982143\sqrt[3]{4}}.
\end{eqnarray}

This next equation is trivial due to it summing a geometric progression, but looks impressive:
\begin{eqnarray}
\sqrt{\frac{5}{9}-\frac{\sqrt[3]{4}}{9}} = \sqrt[3]{\frac{1}{27}} - \sqrt[3]{\frac{2}{27}} + \sqrt[3]{\frac{4}{27}} - \sqrt[3]{\frac{8}{27}} + \sqrt[3]{\frac{16}{27}}.
\end{eqnarray}

%% From "The nested radical identities of Ramanujan"
%% https://www.linkedin.com/grp/post/4510047-5865134162716934149?goback=.gna_4510047
The following nested 10th root equation seems to also be new. It comes easily from the combining of equations listed above.
\begin{eqnarray} \label{E:2.29}
\sqrt[10]{\frac{1}{24} \left(164420\sqrt[3]{25} - 96925\sqrt[3]{\frac{5}{2}} - 273167\sqrt[3]{2}\right)} = \frac{5}{3} \left( \sqrt[3]{\frac{5}{2}} - \sqrt[3]{\frac{1}{5}} + \frac{1}{5} \sqrt[3]{2} \right) \nonumber
\end{eqnarray}

A further new identity is
\begin{eqnarray}
\label{eqx1} && \sqrt[]{4 \sqrt[3]{\frac{2}{3}} - 5 \sqrt[3]{\frac{1}{3}}}, =
\sqrt[3]{\frac{1}{9}}+ 2\sqrt[3]{\frac{2}{9}} - 2\sqrt[3]{\frac{4}{9}}.
\end{eqnarray}
Since we have equal values of various roots of $\sqrt[3]{\frac{1}{9}}-\sqrt[3]{\frac{2}{9}}+\sqrt[3]{\frac{4}{9}}$, these can be used to find further nested examples, of which (\ref{eqx1}) is derived that way from comparing with (\ref{eqx2}).

\section{Constructing nested radicals}

How one comes with such expressions as for example (\ref{eqexa})
\begin{eqnarray*}
 \sqrt[24]{1+100\sqrt[3]{2}-80\sqrt[3]{4}} = \sqrt[3]{\frac{1}{9}}-\sqrt[3]{\frac{2}{9}}+\sqrt[3]{\frac{4}{9}} \; ?
\end{eqnarray*}
Trying to get a root of 24th power from radicals expression is rather difficult.
The "secret" is that we go from RHS of the equation.
We want to express
\begin{eqnarray}
 \label{E:1.2RHS} \sqrt[3]{\frac{1}{9}}-\sqrt[3]{\frac{2}{9}}+\sqrt[3]{\frac{4}{9}}
\end{eqnarray}
as a root of 24th power. So we raise this expression to the 24th power, and after simplification we get
\begin{eqnarray}
  \left(\sqrt[3]{\frac{1}{9}}-\sqrt[3]{\frac{2}{9}}+\sqrt[3]{\frac{4}{9}} \right)^{24} = 1+100\sqrt[3]{2}-80\sqrt[3]{4},
\end{eqnarray}
which is equivalent to (\ref{eqexa}).
\newline

How to make nested radicals equations more "interesting", i.e. having less terms in LHS?

First "method" is random search, or enumeration.
We want to find "interesting" nested radicals with (\ref{E:1.2RHS}) as RHS.
We raise (\ref{E:1.2RHS}) to different powers $n$, starting from $n$ = 2.
We get 2 successful results - 2-term expression - at $n$ = 3, already familiar Ramanujan's equation (\ref{E:1.2}),
and at $n$ = 8:
\begin{eqnarray}
 \sqrt[8]{4\sqrt[3]{\frac{2}{3}}-\frac{5}{\sqrt[3]{3}}} = \sqrt[3]{\frac{1}{9}}-\sqrt[3]{\frac{2}{9}}+\sqrt[3]{\frac{4}{9}}.
\end{eqnarray}

More systematic method. Suppose we want to find an "interesting" nested radical with square root on LHS,
and an expression with $7^{1/4}$, $7^{1/2}$, $7^{3/4}$ and $7$ on RHS.
Simple sum gives
\begin{eqnarray}
  \left(\sqrt[4]{7}+\sqrt[2]{7}+\sqrt[4]{7^3} + 7 \right)^{2} = 28 \sqrt[4]{7} + 22\sqrt{7}+16\sqrt[4]{7^3}+70,
\end{eqnarray}
which is not "interesting". We want on LHS an expression without $\sqrt[4]{7^3}$.
We add coefficients to addends on RHS:
\begin{eqnarray}
&&  \left(\sqrt[4]{7}+x\sqrt[2]{7}+y\sqrt[4]{7^3} + z 7 \right)^{2} =  \\
&&  (14xy+14z) \sqrt[4]{7} + (1+14xz)\sqrt{7}+(2x+14yz)\sqrt[4]{7^3}+(7x^2+14y+49z^2). \nonumber
\end{eqnarray}
We need $(2x + 14yz) = 0$, which may be satisfied at $x$ = 7, $y$ = 1, $z$ = -1, which gives an equation (\ref{E:2.11}).

\end{document}